\documentclass[11pt]{article}
\usepackage{xcolor}
\usepackage{mathrsfs}
\usepackage{amsmath,amsthm,amssymb,amscd}
\usepackage{booktabs} 
\usepackage{float}
\usepackage[hyperfootnotes=false, colorlinks, linkcolor={blue}, citecolor={magenta}, filecolor={blue}, urlcolor={blue}, plainpages=false, pdfpagelabels]{hyperref}
\usepackage[overload]{textcase} 
\usepackage{mathtools}
\usepackage[numbers,sort&compress]{natbib} 
\usepackage{etoolbox} 
\apptocmd{\sloppy}{\hbadness 10000\relax}{}{} 
\usepackage[left=2.5cm, right=2.5cm, top=3cm]{geometry}
\allowdisplaybreaks
\usepackage{enumerate}
\usepackage[shortlabels]{enumitem}
\usepackage{bm}
\usepackage{url}
\usepackage[multiple]{footmisc}

\usepackage{colonequals}
\usepackage{longtable} 
\usepackage[ampersand]{easylist}
\usepackage{ltablex,booktabs}

\def\im{\operatorname{Im}}
\def\re{\operatorname{Re}}

\def\Res{\operatorname{Res}}

\newcommand{\A}{{\mathbb A}}
\newcommand{\Q}{{\mathbb Q}}
\newcommand{\Z}{{\mathbb Z}}

\newcommand{\GL}{{\rm GL}}

\newcommand{\SL}{{\rm SL}}

\newcommand{\Sym}{{\rm Sym}}

\newcommand{\RN}[1]{%
  \textup{\uppercase\expandafter{\romannumeral#1}}%
}

\newcommand{\forget}[1]{}
\def\qdots{\mathinner{\mkern1mu\raise0pt\vbox{\kern7pt\hbox{.}}\mkern2mu
\raise3.4pt\hbox{.}\mkern2mu\raise7pt\hbox{.}\mkern1mu}}

\newtheorem{lemma}{Lemma}[section]
\newtheorem{theorem}[lemma]{Theorem}
\newtheorem{corollary}[lemma]{Corollary}

\newtheorem{remark}[lemma]{Remark}

\newcommand{\pf}{{\noindent\textbf{Proof}\quad}}

\newcommand\blfootnote[1]{%
	\begingroup
	\renewcommand\thefootnote{}\footnote{#1}%
	\addtocounter{footnote}{-1}%
	\endgroup
}
\usepackage[toc,page, title, titletoc]{appendix}
\makeatother
\makeatletter
\newcommand\appendix@section[1]{%
	\refstepcounter{section}%
	\orig@section*{Appendix \@Alph\c@section: #1}%
	\addcontentsline{toc}{section}{Appendix \@Alph\c@section: #1}%
}
\g@addto@macro\appendix{\let\section\appendix@section}
\let\orig@section\section
\makeatother
\makeatletter
\def\thickhline{%
  \noalign{\ifnum0=`}\fi\hrule \@height \thickarrayrulewidth \futurelet
   \reserved@a\@xthickhline}
\def\@xthickhline{\ifx\reserved@a\thickhline
               \vskip\doublerulesep
               \vskip-\thickarrayrulewidth
             \fi
      \ifnum0=`{\fi}}
\makeatother

\title{The Zero Density Theorem for the Rankin-Selberg $L$-function and its applications}
\author{Zhining Wei}
\date{}

\begin{document}

\maketitle
\blfootnote{\hspace*{0.22in} Key words and phrases. Zero density theorem, Rankin-Selberg convolutions, Fourier coefficients of automorphic forms.}

\begin{abstract}
In this work, we establish a zero density result for the Rankin-Selberg $L$-functions. As an application, we apply it to distinguish the holomorphic Hecke eigenforms for $\SL_2(\Z).$
\end{abstract}

\section{Introduction}
The study of zero free region of $L$-functions is a classical theme in analytic number theory. It is widely believed that the generalized Riemann hypothesis (GRH) holds, but a proof for this is out of research at present. In the absence of GRH, the zero density estimates, especially the number of zeros near the line $\re(s)=1$, are often used as a substitute in many applications. Such type results are often referred as the \textit{zero density theorems}. For the zero density theorem of Riemann zeta function $\zeta(s)$, one can refer to \cite[Chapter~12]{Mo1971}. For the zero density theorem for Dirichlet $L$-functions $L(s,\chi),$ one can refer to \cite{Mo1969} and \cite[Theorem~1.4]{IwKo2004}.

Later the zero density theorem for automorphic $L$-functions was well studied in different aspects. In \cite{Lu1999}, Luo established a zero density result for the symmetric square L-function of Maass forms with large eigenvalues. Later Kowalski and Michel \cite{KoMi2002} proved a general density theorem for automorphic $L$-functions with large conductors. In their 2006 work, Lau and Wu \cite{LaWu2006} established a zero density theorem for the holomorphic cusp forms in the weight aspect. Recently Thorner and Zaman \cite{ThZa2021} proved a zero density result for automorphic $L$-functions over a number field $F$ by first establishing an unconditional large sieve inequality for automophic forms on $\GL_n(\A_F).$  For the applications of the zero density theorems, one can also refer to the papers above.

In this paper, we will study the zero density results for the Rankin-Selberg $L$-function of $\SL_2(\Z)$ holomorphic Hecke eigenforms. This will extend the zero free region for Rankin-Selberg functions in the weight aspect \textit{averagely}. As a direct application, we can combine it with \cite{GoHo1993} to distinguish Hecke eigenforms.

We proceed to our results. Let $k$ be an even number. Denote by $\mathcal{S}_k(\SL_2(\Z))$ the (finite dimensional) vector space of weight $k$ holomorphic cusp forms for $\SL_2(\Z).$ Then it has an orthogonal basis with respect to the Petersson inner product, denoted by $H_k,$ which are also eigenfunctions for Hecke operators. Let $f\in H_k$ be a Hecke eigenform. Then it has a Fourier expansion
\[f(z)=\sum_{n=1}^{\infty}\lambda_f(n)n^{\frac{k-1}{2}}e(nz)\]
and we can assume that $\lambda_f(1)=1.$ The associated $L$-function, denoted by $L(s,f)$, has a Euler product:
\[L(s,f)=\sum_{n=1}^{\infty}\frac{\lambda_f(n)}{n^s}=\prod_{p}\left(1-\frac{\alpha_p}{p^s}\right)^{-1}\left(1-\frac{\alpha_p^{-1}}{p^s}\right)^{-1}.\]
Here $\{\alpha_p,\alpha_p^{-1}\}$ are the Satake parameters of $f$ at $p.$ The well-known result by Deligne is $|\alpha_p|=|\alpha_p^{-1}|=1$ for all $p.$ 

Now let $g$ be another Hecke eigenform in $H_k$ with Satake parameters $\{\beta_p,\beta_p^{-1}\}$ for all primes $p.$ Then the Rankin-Selberg $L$-function, denoted by $L(s,f\otimes g),$ is defined by
\[L(s,f\otimes g)=\prod_p\left(1-\frac{\alpha_p\beta_p}{p^s}\right)^{-1}\left(1-\frac{\alpha_p\beta_p^{-1}}{p^s}\right)^{-1}\left(1-\frac{\alpha_p^{-1}\beta_p}{p^s}\right)^{-1}\left(1-\frac{\alpha_p^{-1}\beta_p^{-1}}{p^s}\right)^{-1}\]
when $\re(s)>1.$ The Rankin-Selberg L-function has a meromorphic continuation to the whole complex plane. It can only has a simple pole at $s=1$ and this happens if and only if $f=g.$ In this case, $L(s,f\otimes f)=\zeta(s)L(s,\Sym^2 f)$ with
\[L(s,\Sym^2 f)=\prod_{p}\left(1-\frac{\alpha_p^2}{p^s}\right)^{-1}\left(1-\frac{1}{p^s}\right)^{-1}\left(1-\frac{\alpha_p^{-2}}{p^s}\right)^{-1}\]
when $\re(s)>1.$ Here $L(s,\Sym^2 f)$ is the symmetric square $L$-function assocaited to $f.$ It has an analytic continuation to the whole complex plane.

Next we would introduce a function counting the number of zeros. Let $L(s,\pi)$ be an $L$-function defined in \cite[Chapter~5]{IwKo2004}. ($\pi$ later will be taken to be $\Sym^2 f$ or $f\otimes g.$) For fixed $\alpha>\frac{1}{2}$, denote by $N(\alpha,T,\pi)$ the number of zeros $\rho=\beta+i\gamma$ of $L(s,\pi)$ with $\beta\geq\alpha$ and $0\leq\gamma\leq T.$  We can establish the following zero density results for the Rankin-Selberg $L$-function: 
\begin{theorem}\label{zero density theorem}
Denote by $\mathcal{F}_k$ the set of pairs of distinct Hecke eigenfucntions, that is,
\[\mathcal{F}_k=\{(f,g)\in H_k\times H_k|f\neq g\},\]
For any $\delta>0,$ we have
\[\sum_{(f,g)\in \mathcal{F}_k}N(\alpha,T,f\otimes g)\ll_{\delta} T^2(\log T) k^{34(1-\alpha)/(3-2\alpha)}(\log k)^{25}\]
uniformly for $\frac{1}{2}+\delta\leq \alpha\leq 1$, $(\log k)^{3}\leq T\leq k$ and $f,g\in H_k.$ The constant is dependent on $\delta$ when $\delta$ is small. When $\frac{1}{2}+\delta$ is closed to $1,$ the implied constant is absolute. 
\end{theorem}
In a recent paper \cite{BTZ2022}, Brumley, Thorner and Zaman established a zero density theorem of the Rankin-Selberg $L$-functions under certain Hypothesis. In their paper, they fixed one cuspidal representation $\pi_0$ and varied the other representation in a given family. In our work, we can vary two modular forms (see the definition of $\mathcal{F}_k$). This can be achieved by Lemma \ref{functoriality lemma}. To prove the theorem, we will first establish a large sieve inequality and then follow the work in \cite[Section~5]{LaWu2006}. 

For each $\eta\in (0,1/2),$ define
\[H_k^+(\eta)=\{f\in H_k|\mbox{$L(s,\Sym^2 f)\neq 0$ for $s\in \mathcal{S}$}\}\]
and
\[D^{+}_k(\eta)=\{(f,g)\in\mathcal{F}_k|\mbox{$L(s,f\otimes g)\neq 0$ for $s\in \mathcal{S}$}\},\]
where $S=\{s| \re(s)\geq 1-\eta,|\im(s)\leq 100k^{\eta}|\}\cup\{s|\re(s)\geq 1\}.$

Then set $H_k^-(\eta)=H_k-H_k^+(\eta)$ and $D_k^{-}(\eta)=\mathcal{F}_k-D_k^+(\eta).$ The following corollary shows that $D_k^+(\eta)$ has density one as $k\to\infty$ provided that $\eta$ is small:
\begin{corollary}\label{density one corollary}
For $\eta\in\left(\frac{3\log\log k}{\log k},\frac{1}{4}\right),$ we have
\[|D_k^{-}(\eta)|\ll k^{36\eta}(\log k)^{26}\]
The implied constant is absolute.
\end{corollary}
\noindent\textbf{Proof}: This is similar to \cite[Equation~1.11]{LaWu2006}. When $\eta<\frac{1}{4},$ $\delta\geq\frac{1}{4}.$ The implied constant is absolute due to Remark \ref{Absolute constant}.
\qed
\begin{remark}\label{the modification cor}
A similar argument with Remark \ref{The modification remark} shows that
\[|H_k^-(\eta)|\ll k^{36\eta}(\log k)^{18}.\]
\end{remark}

Next we will consider an application of the zero density theorem: let $f,g\in H_k$ be distinct Hecke eigenforms. We say that $f$ and $g$ are \textit{distinguishable} if for every $\epsilon>0,$ we can find $n\ll_{\epsilon} k^{\epsilon}$ such that $\lambda_f(n)\neq\lambda_g(n)$. 
Combine Corollary \ref{density one corollary} and Remark \ref{the modification cor}, and we establish the following result:
\begin{theorem}\label{distinguish theorem}
As $k\to\infty,$ we can find a set of Hecke eigenforms $H_k^*$ ($\subseteq H_k$) such that
\[\lim_{k\to\infty}\frac{|H_k^*|}{|H_k|}=1\]
and for any $f,g\in H_k^*$ and $f\neq g,$ they are distinguishable. 
\end{theorem}
The proof is based on \cite{GoHo1993}. 

\section{The large sieve inequality}
In this section, we will prove the following large sieve inequality, which will be used to prove Theorem \ref{zero density theorem} in the next section. Before that, we need the following lemma to know the possible poles for Rankin-Selberg $L$-functions:
\begin{lemma}\label{functoriality lemma}
Let $f_1,f_2,g_1,g_2\in H_k$ be Hecke eigenforms. Assume that $f_1\neq g_1$ and $f_2\neq g_2.$ Then $L(s,f_1\otimes g_1\otimes f_2\otimes g_2)$ has at most a simple pole at $s=1.$ This happens if and only one of the following is valid:
\begin{enumerate}[(a)]
\item 
$f_1=f_2$ and $g_1=g_2$;
\item
$f_1=g_2$ and $f_2=g_1.$
\end{enumerate}
\end{lemma}
\begin{proof}
($\Rightarrow$) By the condition that $f_1\neq g_1$ and $f_2\neq g_2,$ we know that the $L(s,f_1\otimes g_1\otimes f_2\otimes g_2)$ will have a pole at $s=1$ if one of (a), (b) holds. 

($\Leftarrow$) If exactly two of $f_1,f_2,g_1,g_2$ are equal, then $L(s,f_1\otimes g_1\otimes f_2\otimes g_2)$ is entire. Without loss of generality, we assume that $f_1=f_2$ and $f_1,g_1,g_2$ are distinct. Then 
\[L(s,f_1\otimes g_1\otimes f_2\otimes g_2)=L(g_1\otimes g_2)L(s,\Sym^2 f_1\otimes g_1\otimes g_2)\]
is entire due to \cite{Ra2000} and \cite[Theorem~9.2]{CoKiMu2004}. The left cases can be discussed similarly. So it suffices to consider the case when any two of $f_1,f_2,g_1,g_2$ are distinct. Then by \cite{Ra2000}, we know that, for each pair of distinct $f,g\in\{f_1,f_2,g_1,g_2\},$ there exists an automorphic cuspidal representations of $\GL_4(\A_{\Q}),$ denoted by $\Pi_{f,g}$, such that $L(s,f\otimes g)=L(s,\Pi_{f,g})$. Then by \cite[Theorem~9.2]{CoKiMu2004}, $L(s,f_1\otimes g_1\otimes f_2\otimes g_2)$ will have at most a simple pole at $s=1.$ We prove the rest by contradiction: suppose that $L(s,f_1\otimes g_1\otimes f_2\otimes g_2)$ has a simple pole at $s=1.$ Then for any $x\geq 1,$ we can show that
\begin{equation}\label{simple pole}
    \sum_{p\leq x}\lambda_{f_1}(p)\lambda_{f_2}(p)\lambda_{g_1}(p)\lambda_{g_2}(p)=\frac{x}{\log x}+O_A\left(\frac{x}{\log^A x}\right)
\end{equation}
by \cite[Theorem~5.13]{IwKo2004}. Here $A>0$ is arbitrary. (Notice that $f_1,f_2,g_1,g_2$ satisfy the Ranmanujan conjecture and hence \cite[Equation~(5.48)]{IwKo2004} is obvious.) 

On the other hand, \cite[Theorem~9.2]{CoKiMu2004} implies that
\[
\Pi_{f_1,f_2}\cong \Pi_{g_1,g_2}\hspace{6mm} \Pi_{f_1,g_1}\cong \Pi_{f_2,g_2} \hspace{6mm} \Pi_{f_1,g_2}\cong \Pi_{f_2,g_1}.
\]
Therefore, for each prime $p,$ we have
\[\lambda_{f_1}(p)\lambda_{f_2}(p)=\lambda_{g_1}(p)\lambda_{g_2}(p) \hspace{6mm} \lambda_{f_1}(p)\lambda_{g_1}(p)=\lambda_{f_2}(p)\lambda_{g_2}(p) \hspace{6mm} \lambda_{f_1}(p)\lambda_{g_2}(p)=\lambda_{f_2}(p)\lambda_{g_1}(p).\]
Multiply them together, and we obtain that
\[\lambda_{f_1}(p)^3\lambda_{f_2}(p)\lambda_{g_1}(p)\lambda_{g_2}(p)=\lambda_{f_2}(p)^2\lambda_{g_1}(p)^2\lambda_{g_2}(p)^2.\]
Let $S$ be a set of primes satisfying
\[S=\{p|\lambda_{f_2}(p)\lambda_{g_1}(p)\lambda_{g_2}(p)=0\}.\]
Since $L(s,f_1\otimes f_1\otimes f_1\otimes f_1)$ has a pole of order $2$ at $s=1,$ we have
\begin{equation}\label{secodn order pole}
\begin{split}
\sum_{p\leq x}\lambda_{f_1}(p)\lambda_{f_2}(p)\lambda_{g_1}(p)\lambda_{g_2}&=\sum_{p\leq x,p\notin S}\lambda_{f_1}(p)\lambda_{f_2}(p)\lambda_{g_1}(p)\lambda_{g_2}(p)\\
&\geq\sum_{p\leq x}\lambda_{f_1}(p)^4-16\sum_{p\leq x,p\in S}1\hspace{8mm}(|\lambda_{f_1}(p)|\leq2)\\
&=\frac{2x}{\log x} +O_A\left(\frac{x}{\log^A x}\right)-16\sum_{p\leq x,p\in S}1
\end{split}
\end{equation}
However $p\in S$ implies that one of $\lambda_{f_2}(p),\lambda_{g_1}(p),\lambda_{g_2}(p)$ is zero. This cannot happen quite frequently due to the Sato-Tate conjecture \cite{BGHT2011}. ($\lambda_{f}(p)=0$ implies that one of Satake parameters $\alpha_p,\alpha_p^{-1}$ is $e^{i\pi/2}$) Indeed, we can show that for any $\epsilon>0,$
\[\sum_{p\leq x, p\in S}1\leq\frac{\epsilon x}{\log x}\]
as $x\to\infty.$ Then combine this with Equation \eqref{simple pole} and \eqref{secodn order pole}, a contradiction.
\end{proof}

Then we prove the following large sieve inequality for the family $\mathcal{F}_k:$
\begin{lemma}\label{large sieve inequality}
Let $L\geq 1$. Let $\{a_{\ell}\}_{\ell\leq L}$ be a sequence of complex numbers. Then for any $\epsilon>0,$ we have
\[\sum_{(f,g)\in\mathcal{F}_k}\left|\sum_{\ell\leq L}a_{\ell}\lambda_{f\otimes g}(\ell)\right|^2\ll_{\epsilon}(L(\log k)^{15}+k^{9/2+\epsilon}L^{1/2+\epsilon})\sum_{\ell\leq L}|a_{\ell}|^2.\]

\end{lemma}
\pf By duality principal, it suffices to show:
\[\sum_{\ell\leq L}\left|\sum_{(f,g)\in \mathcal{F}_k}b_{f,g}\lambda_{f\otimes g}(\ell)\right|^2\ll_{\epsilon}(L\log L(\log k)^{15}+k^{9/2+\epsilon}L^{1/2+\epsilon})\sum_{f,g\in H_k}|b_{f,g}|^2.\]

The left hand side is
\begin{equation}\label{duality sum}
\ll\sum_{\ell\geq 1}\left|\sum_{(f,g)\in \mathcal{F}_k}b_{f,g}\lambda_{f\otimes g}(\ell)\right|^2e^{-\ell/L}=\sum_{\substack{(f_1,g_1)\in \mathcal{F}_k\\(f_2,g_2)\in \mathcal{F}_k}}b_{f_1,g_1}\overline{b_{f_2,g_2}}\sum_{\ell\geq 1}\lambda_{f_1\otimes g_1}(\ell)\lambda_{f_2\otimes g_2}(\ell)e^{-\ell/L}.
\end{equation}
A standard argument will show that
\begin{flalign*}
\sum_{\ell\geq 1}\lambda_{f_1\otimes g_1}(\ell)\lambda_{f_2\otimes g_2}(\ell)e^{-\ell/L}&=\frac{1}{2\pi i}\int_{(2)}L(s,f_1\otimes g_1\otimes f_2\otimes g_2)G_{f_1,f_2,g_1,g_2}(s)\Gamma(s)L^s\,ds\\
&=\Res\limits_{s=1}L(s,f_1\otimes g_1\otimes f_2\otimes g_2)G_{f_1,f_2,g_1,g_2}(s)\Gamma(s)L^s\\
&\hspace{6mm}+\frac{1}{2\pi i}\int_{(1/2+\epsilon)}L(s,f_1\otimes g_1\otimes f_2\otimes g_2)G_{f_1,f_2,g_1,g_2}(s)\Gamma(s)L^s\,ds
\end{flalign*}
Here $G_{f_1,f_2,g_1,g_2}(s)$ is some Euler product which is absolutely convergent for $\re(s)>1/2.$ Additionally, we have $G_{f_1,f_2,g_1,g_2}(s)\ll_{\epsilon} 1$ for $\re(s)\geq\frac{1}{2}+\epsilon$ (independent from the choice of $f_1,f_2,g_1$ and $g_2$ since they satisfy the Ranmanujan conjecture.) By the virtual of the proof bellow, it suffices to assume that $G_{f_1,f_2,g_1,g_2}(1)\neq 0.$

For the function $L(s,f_1\otimes g_1\otimes f_2\otimes g_2),$ it will have a pole of order $1$ if and only if either $f_1=f_2\neq g_1=g_2$ or $f_1=g_2\neq f_2=g_1$ by Lemma \ref{functoriality lemma}.  In both cases, we have
\[L(s,f_1\otimes g_1\otimes f_2\otimes g_2)=\zeta(s)L(s,\Sym^2f_1)L(s,\Sym^2g_1)L(s,\Sym^2f_1\otimes \Sym^2g_1).\]

Since $f_1,f_2,g_1,g_2$ satisfies the Ranmanujan conjecture, we can show:
\[L(1,\Sym^2f)\ll(\log k)^{3}\]
and
\[L(1,\Sym^2f\otimes\Sym^2g)\ll(\log k)^{9}\]
by \cite[Lemma~4.1]{CoMi2004}. 
This will show that
\[\Res\limits_{s=1}L(s,f_1\otimes g_1\otimes f_2\otimes g_2)G_{f_1,f_2,g_1,g_2}(s)\Gamma(s)L^s\ll (\delta_{f_1,f_2}\delta_{g_1,g_2}+\delta_{f_1,g_2}\delta_{f_2,g_1})L(\log k)^{15}.\]
Here $\delta_{f,g}=1$ if $f=g$ and $0$ otherwise. By \cite[Section~3]{LaWu2006}, one can show that
\begin{equation}\label{convexity bound}
L(1/2+\epsilon+it,f_1\otimes g_1\otimes f_2\otimes g_2)\ll_{\epsilon}(1+|t|)^{6/4}(k+|t|)^{5/2+\epsilon}.
\end{equation}
This will show that
\begin{equation}\label{inner term for large sieve}
\sum_{\ell\geq 1}\lambda_{f\otimes g_1}(\ell)\lambda_{f\otimes g_2}(\ell)e^{-\ell/L}\ll_{\epsilon}(\delta_{f_1,f_2}\delta_{g_1,g_2}+\delta_{f_1,g_2}\delta_{f_2,g_1})L(\log k)^{15}+k^{5/2+\epsilon}L^{1/2+\epsilon}.
\end{equation}
Notice that, by Cauchy's inequality,
\[\sum_{(f,g)\in\mathcal{F}_k}|b_{f,g}\overline{b_{g,f}}|\leq \left(\sum_{(f,g)\in\mathcal{F}_k}|b_{f,g}|^2\right).\]
Insert Equation \eqref{inner term for large sieve} to Equation \eqref{duality sum} and we obtain the result.
\qed
\begin{remark}\label{modification of the main term}
By a similar argument, one can show that the term $Lk^{\epsilon}$ in \cite[Proposition~4.1]{LaWu2006} can be replaced by $L(\log k)^8.$
\end{remark}

\section{Proof of Theorem \ref{zero density theorem}}
Our proof is based on the method of Montgomery in \cite{Mo1971}. Here we will follow the work in \cite[Section~5]{LaWu2006}. We will modify their work to replace $k^{\epsilon}$ by a large power of $\log k.$ Indeed, this also works for the symmetric square $L$-functions and we will discuss it in Remark \ref{The modification remark}.

We first establish the following lemma, which is similar to \cite[Lemma~5.1]{LaWu2006}:
\begin{lemma}\label{approximation lemma}
Let $z>16$ be any fixed number and let $P(z)=\prod_{p\leq z}p.$ For any $\re(s)=\sigma>1,$ we have
\[L(s,f\otimes g)^{-1}=G_{f,g}(s)\sum_{(n,P(z))=1}\frac{\lambda_{f\otimes g}(n)\mu(n)}{n^s}.\]
The Dirichlet series $G_{f,g}(s)$ converges absolutely for $\sigma>\frac{1}{2}$ and $G_{f,g}(s)\ll_{\epsilon} 1$ uniformly for $\re(s)\geq\frac{1}{2}+\epsilon.$ Notice that the implied constant is independent from the choice of $f,g.$
\end{lemma}

Then for $\re(s)>1/2,$ we define
\[M_x(s,f\otimes g)=G_{f,g}(s)\sum_{\ell\leq x,(\ell,P(z))=1}\frac{\mu(\ell)\lambda_{f\otimes g}(\ell)}{\ell^s}.\]
The trivial bound shows that, for any $\epsilon>0$ 
\[M_x(s,f\otimes g)\ll_{\epsilon} x^{1/2}\]
uniformly for $\re(s)\geq\frac{1}{2}+\epsilon.$ Obviously, we have
\[1=(1-L(s,f\otimes g)M_x(s,f\otimes g))+L(s,f\otimes g)M_x(s,f\otimes g).\]
\vspace{6mm}

\noindent\textbf{Proof of Theorem \ref{zero density theorem}}: 
We cut the rectangle $\alpha\leq \re(s)\leq 1$ and $0\leq \im(s)\leq T$ horizontally into boxes of width $2(\log k)^2.$ Then by \cite[Propositon~5.7]{IwKo2004}, each box $\alpha\leq \re(s)\leq 1$ and $Y\leq \im(s)\leq Y+2(\log k)^2$ contains at most $O((\log k)^3)$ zeros. Denote by $n_{f\otimes g}$ the number of boxes which contains at least a zero of $L(s,f\otimes g).$ Then
\[N(\alpha, T,f\otimes g)\ll n_{f\otimes g}(\log k)^3\ll n_{f\otimes g}T.\]
So to prove Theorem \ref{zero density theorem}, it suffices to show: for a fixed $\alpha\geq \frac{1}{2}+\delta,$
\[\sum_{(f,g)\in \mathcal{F}_k}n_{f\otimes g}\ll_{\delta} T(\log T)k^{34(1-\alpha)/(3-2\alpha)}(\log k)^{25}.\]

Let $x,y$ be large numbers to be chosen later. (We will choose $x,y$ to be some power of $k$ and such choice is sufficient for the following argument.) 
Let $\rho=\beta+i\gamma$ with $\beta\geq\alpha(>\frac{1}{2}+\epsilon)$ and $0\leq \gamma\leq T$ be a zero of $L(s,f\otimes g)$ and we write:
\[\kappa=\frac{1}{\log k}\hspace{10mm}\kappa_1=1-\beta+\kappa(>0)\hspace{10mm}\kappa_2=\frac{1}{2}-\beta+\epsilon(<0).\]
(Here it suffices to assume that $\epsilon<\delta.$) Then follow \cite[Section~5]{LaWu2006}, we have:
\begin{flalign*}
e^{-1/y}&=\frac{1}{2\pi i}\int_{(\kappa_1)}(1-L(\rho+\omega,f\otimes g)M_x(\rho+\omega,f\otimes g))\Gamma(\omega)y^{\omega}\,d\omega\\
&\hspace{6mm}+\frac{1}{2\pi i}\int_{(\kappa_2)}L(\rho+\omega,f\otimes g)M_x(\rho+\omega,f\otimes g)\Gamma(\omega)y^{\omega}\,d\omega\\
\end{flalign*}
By the convexity bound in \cite[Proposition~3.1]{LaWu2006}: 
\begin{equation}\label{convexity bound for Rankin-Selberg}
L(s,f\otimes g)\ll_{\epsilon'}(1+|\im(s)|)^{1-\re(s)}(k+|\im(s)|)^{1-\re(s)+\epsilon'}
\end{equation}
for any $\epsilon'>0,$ we have:
\[\int\limits_{(\kappa_1), |\im(\omega)|\geq (\log k)^2}(1-L(\rho+\omega,f\otimes g)M_x(\rho+\omega,f\otimes g))\Gamma(\omega)y^{\omega}\,d\omega\ll_{\epsilon} \frac{1}{k^2}\]
and
\[\int\limits_{(\kappa_2), |\im(\omega)|\geq (\log k)^2}L(\rho+\omega,f\otimes g)M_x(\rho+\omega,f\otimes g)\Gamma(\omega)y^{\omega}\,d\omega\ll_{\epsilon}\frac{1}{k^2}.\]
This is due to the rapid decay of $\Gamma(\omega)$ when $|\im(\omega)|\geq(\log k)^2.$ In this case, set $K=(\log k)^2$ and by the fact that $1\leq C(a+b)\Rightarrow 1\leq 2C^2(a^2+b)$ (where $a,b>0$ and $C\geq 1,$), we obtain:
\begin{flalign*}
1&\ll (\log k)^2y^{2-2\alpha}\int_{-K}^{K}|1-L(1+\kappa+i(\gamma+v),f\otimes g)M_x(1+\kappa+i(\gamma+v),f\otimes g)|^2\,dv\\
&\hspace{6mm}+y^{1/2-\alpha+\epsilon}\int_{-K}^K|L(1/2+\epsilon+i(\gamma+v),f\otimes g)M_x(1/2+\epsilon+i(\gamma+v),f\otimes g)|\,dv.
\end{flalign*}
Then follow the work in \cite[Page~457]{LaWu2006}, and we have:
\begin{flalign*}
n_{f\otimes g}&\ll (\log k)^2y^{2-2\alpha}\int_{0}^{2T}|1-L(1+\kappa+iv,f\otimes g)M_x(1+\kappa+iv,f\otimes g)|^2\,dv\\
&\hspace{6mm}+y^{1/2-\alpha+\epsilon}\int_{0}^{2T}|L(1/2+\epsilon+iv,f\otimes g)M_x(1/2+\epsilon+iv,f\otimes g)|\,dv\\
&=:(\log k)^2y^{2-2\alpha} \RN{1}_{f\otimes g}+y^{1/2-\alpha+\epsilon}\RN{2}_{f\otimes g}.
\end{flalign*}
By Equation \eqref{convexity bound for Rankin-Selberg}, we have, for any $\epsilon'>0,$
\[L(1/2+\epsilon+iv,f\otimes g)\ll_{\epsilon'}(1+|v|)^{1/2}(k+|v|)^{1/2+\epsilon'}\ll_{\epsilon'} k^{1+\epsilon'}\]
provided that $T\leq k.$ Therefore, for $\RN{2}_{f\otimes g},$ we have:
\begin{equation}\label{II estimate}
\RN{2}_{f\otimes g}\ll_{\epsilon,\epsilon'} Tx^{1/2}k^{1+\epsilon'}.
\end{equation}

Then set $X=e^{4(\log k)^2}$ and we have
\begin{flalign*}
1-&L(1+\kappa+iv,f\otimes g)M_x(1+\kappa+iv,f\otimes g)\\
&=L(1+\kappa+iv,f\otimes g)G_{f,g}(1+\kappa+iv)\sum_{\ell>x,(\ell,P(z))=1}\frac{\mu(\ell)\lambda_{f\otimes g}(\ell)}{\ell^{1+\kappa+iv}}\\
&\ll_{\epsilon} \zeta(1+\kappa)^4\left(\left|\sum_{x<\ell\leq X,(\ell,P(z))=1}\frac{\mu(\ell)\lambda_{f\otimes g}(\ell)}{\ell^{1+\kappa+iv}}\right|+\sum_{\ell>X}\frac{d_4(\ell)}{\ell^{1+\kappa}}\right)\\
&\ll_{\epsilon}\zeta(1+\kappa)^4\left|\sum_{x<\ell\leq X,(\ell,P(z))=1}\frac{\mu(\ell)\lambda_{f\otimes g}(\ell)}{\ell^{1+\kappa+iv}}\right|+X^{-\kappa/2}\zeta(1+\kappa/2)^4\zeta(1+\kappa)^4
\end{flalign*}
Recall that $\kappa=\frac{1}{\log k}$ and this will give:
\[\sum_{(f,g)\in \mathcal{F}_k}\RN{1}_{f\otimes g}\ll_{\epsilon}(\log k)^4\int_0^{2T}\sum_{(f,g)\in \mathcal{F}_k}\left|\sum_{x<\ell\leq X,(\ell,P(z))=1}\frac{\mu(\ell)\lambda_{f\otimes g}(\ell)}{\ell^{1+\kappa+iv}}\right|^2\,dv+T(\log k)^8.\]
Then by dyadic division, Lemma \ref{large sieve inequality} and Cauchy-Schwarz inequality, we obtain:
\begin{equation}\label{I estimate}
\sum_{(f,g)\in \mathcal{F}_k}\RN{1}_{f\otimes g}\ll_{\epsilon,\epsilon_1}(\log k)^4(\log X)^2T\left((\log k)^{15}+k^{9/2+\epsilon_1}x^{-1/2+\epsilon_1}\right).
\end{equation}
Combine Equation \eqref{I estimate} and Equation \eqref{II estimate}, and we have:
\begin{flalign*}
\sum_{(f,g)\in \mathcal{F}_k}n_{f\otimes g}&\ll_{\epsilon}(\log k)^2y^{2-2\alpha} \sum_{(f,g)\in \mathcal{F}_k}\RN{1}_{f\otimes g}+y^{1/2-\alpha+\epsilon}\sum_{(f,g)\in \mathcal{F}_k}\RN{2}_{f\otimes g}\\
&\ll_{\epsilon,\epsilon_1,\epsilon'} (\log k)^2y^{2-2\alpha}\left(T(\log k)^{23}+Tk^{9/2+\epsilon_1}x^{-1/2+\epsilon_1}(\log k)^{4}\right)+y^{1/2-\alpha+\epsilon}Tx^{1/2}k^{3+\epsilon'}
\end{flalign*}
Set $\epsilon_1=\frac{1}{22}, \epsilon'=\frac{1}{4}$ and $\epsilon=\min\left\{\frac{\delta}{2},\frac{1}{100}\right\}$. Then take $x=k^{10}$
and $y=k^{17/(3-2\alpha)}$ and we obtain the result.
\qed

\begin{remark}\label{Absolute constant}
It can be seen from the proof that, the implied constant comes from the function $G_{f,g}(s)$. When we consider $\delta\geq\frac{1}{4},$ $\epsilon=\frac{1}{100}.$ So the implied constant is absolute.
\end{remark}

\begin{remark}\label{The modification remark}
For the symmetric square $L$-function $L(s,\Sym^2 f),$ one can show that for $\alpha\geq\frac{1}{2}+\delta,$
\[\sum_{f\in H_k}N(\alpha,T,\Sym^2f)\ll_{\delta}T^{2}(\log T) k^{22(1-\alpha)/(3-2\alpha)}(\log k)^{17}.\]
Furthermore, when $\frac{1}{2}+\delta$ is close to $1,$ the implied constant is absolute.
\end{remark}

\section{Proof of Theorem \ref{distinguish theorem}}
\noindent\textbf{Proof of Theorem \ref{distinguish theorem}}: Fix $\epsilon>0.$ Set
\[\eta=\frac{(\log\log k)^2}{\log k}\in\left(\frac{3\log\log k}{\log k},\frac{1}{4}\right),\]
(This is true when $k$ is large.) and we apply the Remark \ref{the modification cor} and Corollary \ref{density one corollary}. This gives:
\begin{equation}\label{symmetric square exceptions}
|H_k^-(\eta)|\ll k^{36\eta}(\log k)^{18}=o((\log k)^{37\log\log k})
\end{equation}
and 
\begin{equation}\label{Rankin-Selberg exceptions}
|D_k^-(\eta)|\ll k^{36\eta}(\log k)^{26}=o((\log k)^{37\log\log k}).
\end{equation}

Then we claim: for $f\in H_k^+(\eta)$ and $(f,g)\in D_k^+(\eta)$, $f$ and $g$ are distinguishable. We will prove the claim later. Then we define 
\[H_k(-,\eta)=\{f\in H_k|\mbox{there exists $g\in H_k$ such that either $(f,g)$ or $(g,f)$ belongs to $D_k^-(\eta)$}\}\]
Then equation \eqref{Rankin-Selberg exceptions} implies that
\begin{equation}\label{Bad term}
|H_k(-,\eta)|\leq |D_k^-(\eta)|=o((\log k)^{37\log\log k}).
\end{equation}
Then set
\[H_k^*=H_k-(H_k^-(\eta)\cup H_k(-,\eta))\]
It can be seen that for distinct $f,g\in H_k^*,$ we have $f\in H_k^+(\eta)$ and $(f,g)\in D_k^+(\eta).$ If we assume the claim, then $f$ and $g$ are distinguishable. Then by Equation \eqref{symmetric square exceptions} and \eqref{Bad term},  we have
\[\lim_{k\to\infty}\frac{|H_k^*|}{|H_k|}=1.\]

So it suffices to prove the claim. The main idea is to compare two integrals, which comes from \cite{GoHo1993}. Let $f\in H^+(\eta)$ and $(f,g)\in D_k^+(\eta),$ $g\neq f.$ For $\re(s)>1,$ we define $\Lambda_{f\otimes f}(n)$ and $\Lambda_{f\otimes g}(n)$ by:
\[-\frac{L'(s,f\otimes f)}{L(s,f\otimes f)}=\sum_{n=1}^{\infty}\frac{\Lambda_{f\otimes f}(n)}{n^s}\]
and
\[-\frac{L'(s,f\otimes g)}{L(s,f\otimes g)}=\sum_{n=1}^{\infty}\frac{\Lambda_{f\otimes g}(n)}{n^s}.\]
Then set $a=1+\frac{1}{\log k}$ and $T=50k^{\eta}=50(\log k)^{\log\log k}$. For any $x\geq 1$, we consider
\[\RN{1}=\frac{1}{2\pi i}\int_{a-iT}^{a+iT}\left(\frac{x^{s-1/2}-x^{1/2-s}}{s-1/2}\right)^2\left(-\frac{L'(s,f\otimes f)}{L(s,f\otimes f)}\right)\,ds.\]
This gives:
\begin{equation}\label{RN1 left}
    \RN{1}=\sum_{n<x^2}\frac{\Lambda_{f\otimes f}(n)}{n^{1/2}}\log\left(\frac{x^2}{n}\right)+O\left(\frac{x^{2a-1}(\log k)}{T}\right)
\end{equation}
since
\[\frac{1}{2\pi i}\int_{a-i\infty}^{a+i\infty}\frac{y^{s-1/2}}{(s-1/2)^2}\,ds=\left\{\begin{array}{cc}
\log y&\mbox{if $y\geq 1$}\\
0&\mbox{if $y\leq 1$}.
\end{array}
\right.\]
On the other hand, we consider the path which is the boundary of the rectangular:
\[\left\{z\in\mathbb{C}\left|1-\frac{3}{4}\eta\leq\re(s)\leq a, |\im(s)|\leq T\right.\right\}.\]
Due to the choice of $f$, $L(s,f\otimes f)=\zeta(s)L(s,\Sym^2 f)$ has no zero in the rectangular: the rectangle is contained in the zero free region of Riemann zeta function when $k$ large. Here we recall the well--known zero free region of Riemann zeta function proved by Vinogradov-Korobov \cite{Vi1958} is, for $s=\beta+it$ and
\[\beta\in\left(1-\frac{c}{(\log t)^{2/3}(\log\log t)^{1/3}},1\right),\]
$\zeta(s)\neq 0.$ On the other hand, since $f\in H_k(\eta),$ $L(s,\Sym^2 f)$ has no zeros in the rectangle. (Notice that the $L$-function is self-dual.) There is a simple pole of $L(s,f\otimes f)=\zeta(s)L(s,\Sym^2 f)$ at $s=1.$ Then by the residue theorem, we have
\begin{flalign*}
\RN{1}&=4(x-2+x^{-1})\\
&\hspace{6mm}+\frac{1}{2\pi i}\left(\int_{1-\frac{3}{4}\eta-i\infty}^{1-\frac{3}{4}\eta+i\infty}+\int_{1-\frac{3}{4}\eta+iT}^{a+iT}-\int_{1-\frac{3}{4}\eta-iT}^{a-iT}\right)\left(\frac{x^{s-1/2}-x^{1/2-s}}{s-1/2}\right)^2\left(-\frac{L'(s,f\otimes f)}{L(s,f\otimes f)}\right)\,ds\\
&=4(x-2+x^{-1})+\RN{1}_1+\RN{1}_2+\RN{1}_3.
\end{flalign*}
Then by \cite[Proposition~5.7]{IwKo2004}, we have:
\begin{equation}\label{RN1 right23}
    \RN{1}_2+\RN{1}_3\ll\frac{x^{2a-1}(\log k)^2}{T^2}.
\end{equation}
For the term $\RN{1}_1,$ we have the trivial estimation:
\begin{equation}\label{RN1 right1}
    \RN{1}_1\ll x^{1-\frac{3}{2}\eta}(\log k)^2
\end{equation}
Then combine Equation \eqref{RN1 left}, \eqref{RN1 right23} and \eqref{RN1 right1}, and we obtain:
\begin{equation}\label{RN1 estimate}
    \sum_{n<x^2}\frac{\Lambda_{f\otimes f}(n)}{n^{1/2}}\log\left(\frac{x^2}{n}\right)=4(x-2+x^{-1})+O\left(\frac{x^{2a-1}(\log k)^2}{T}\right)+O(x^{1-\frac{3}{2}\eta}(\log k)^2).
\end{equation}
Then we replace $-L'(s,f\otimes f)/L(s,f\otimes f)$ by $-L'(s,f\otimes g)/L(s,f\otimes g)$, that is, set
\[\RN{2}=\frac{1}{2\pi i}\int_{a-iT}^{a+iT}\left(\frac{x^{s-1/2}-x^{1/2-s}}{s-1/2}\right)^2\left(-\frac{L'(s,f\otimes g)}{L(s,f\otimes g)}\right)\,ds.\]
This will give:
\begin{equation}\label{RN2 estimate}
    \sum_{n<x^2}\frac{\Lambda_{f\otimes g}(n)}{n^{1/2}}\log\left(\frac{x^2}{n}\right)=O\left(\frac{x^{2a-1}(\log k)^2}{T}\right)+O(x^{1-\frac{3}{2}\eta}(\log k)^2).
\end{equation}
since $L(s,f\otimes g)$ has neither zero nor poles in the rectangular.

Now suppose that $\Lambda_{f\otimes f}(n)=\Lambda_{f\otimes g}(n)$ for $n<x^2.$ Then Equation \eqref{RN1 estimate} and \eqref{RN2 estimate} show:
\[0=4(x-2+x^{-1})+O\left(\frac{x^{2a-1}(\log k)^2}{T}\right)+O(x^{1-\frac{3}{2}\eta}(\log k)^2).\]
For any $\epsilon>0$, we set $x=k^{\epsilon/2}.$ By the choice of $c,\eta$ and $T,$ this can not be true when $k$ is large. So we can find $n\leq x^2=k^{\epsilon}$ such that $\Lambda_{f\otimes f}(n)\neq\Lambda_{f\otimes g}(n).$ Notice that $\Lambda_{f\otimes f}(n)$ and $\Lambda_{f\otimes g}(n)$ are supported on prime powers and totally determined by Satake parameters. Therefore, we can find $p\ll_{\epsilon} k^{\epsilon}$ such that $\lambda_f(p)\neq\lambda_g(p)$.
\qed

\section*{Acknowledgments}
The author would like to thank Professor Wenzhi Luo for the suggestions on the topic.

\newpage 

\bibliographystyle{alpha}
\bibliography{ZDA.bib}
\end{document}